\documentclass[a4paper,11pt, fleqn]{amsart}

\usepackage{pst-plot}

\usepackage{amssymb}
\usepackage{amsmath}
\usepackage{enumerate}

\usepackage[all]{xy}
\theoremstyle{plain}
\newtheorem{theorem}{Theorem}

\newtheorem{proposition}{Proposition}
\newtheorem{lemma}[proposition]{Lemma}

\theoremstyle{definition}
\newtheorem{definition}[proposition]{Definition}
\newtheorem{example}[proposition]{Example}

\theoremstyle{remark}
\newtheorem{remark}[proposition]{Remark}

\newcommand{\secref}[1]{Section~\ref{#1}}

\newcommand{\thmref}[1]{Theorem~\ref{#1}}
\newcommand{\propref}[1]{Proposition~\ref{#1}}

\newcommand{\defref}[1]{Definition~\ref{#1}}

\def\cA{{\mathcal A}}

\def\cF{{\mathcal F}}

\def\cS{{\mathcal S}}

\def\R{\mathbb{R}}

\def\id{{\rm id}}
\def\proj{{\rm proj}}
\newcommand{\Top}{\mbox{\bf Top}}
\newcommand{\Tops}{\mbox{${\mathcal S}$-$\bf Top$}}

\newcommand{\Gcat}{{\rm Gcat}\,}
\newcommand{\Ocat}{{\rm Ocat}\,}
\newcommand{\catF}{{\rm cat_\cF}\,}

\newcommand{\cat}{{\rm cat}\,}

\newcommand{\Whcat}{{\rm Whcat}\,}
\newcommand{\SO}{{\rm SO}\,}
\newcommand{\npu}{{n+1}}

\newcommand{\transverse}{\pitchfork}
\renewcommand{\smallskip}{\vskip2mm}
\renewcommand{\medskip}{\vskip4mm}
\begin{document}
\title[Ganea and Whitehead definitions for tangential LS-category]{Ganea and Whitehead definitions for the tangential Lusternik-Schnirelmann category of foliations
}

\author[J.P Doeraene]{Jean-Paul Doeraene}
\address{D\'epartement de Math\'ematiques\\
        UMR-CNRS 8524\\
        Universit\'e de Lille~1\\
        59655 Villeneuve d'Ascq Cedex\\
        France}
\email{Jean-Paul.Doeraene@univ-lille1.fr}

\author[E. Macias-Virg{\'o}s]{Enrique Macias-Virg{\'o}s}
\address{Departamento de Xeometria e Topoloxia\\
        Facultade de Matem{\'a}ticas\\
        Universidade de Santiago de Compostela\\
			      15706 - Galicia\\
        Spain}
\email{macias@zmat.usc.es}
\thanks{The second author is partially supported by research project MICINN MTM2008-05861}

\author[D. Tanr{\'e}]{Daniel Tanr{\'e}}
\address{D\'epartement de Math\'ematiques\\
        UMR-CNRS 8524\\
        Universit\'e de Lille~1\\
        59655 Villeneuve d'Ascq Cedex\\
        France}
\email{Daniel.Tanre@univ-lille1.fr}

\begin{abstract}
This work solves the problem of elaborating Ganea and Whitehead definitions for the tangential category of a foliated manifold.  We develop  these two notions in  the category $\Tops$ of stratified spaces, that are topological spaces $X$ endowed with a partition $\cF$ and compare them to a third invariant defined by using open sets.  More precisely, these definitions apply to an element $(X,\cF)$ of $\Tops$ together with a class $\cA$ of subsets  of $X$; they are similar to invariants introduced by M.~Clapp and D.~Puppe.

If $(X,\cF)\in\Tops$, we define a transverse subset   as a subspace $A$ of $X$ such that the intersection $S\cap A$ is at most countable for any $S\in \cF$.  
Then we define the Whitehead and Ganea LS-categories of the stratified space by taking the infimum along the transverse subsets.
When we have a closed manifold, endowed with a $C^1$-foliation, the three previous definitions, with $\cA$ the class of transverse subsets, coincide with the tangential category and are homotopical invariants.
\end{abstract}



\maketitle

\tableofcontents

\newpage
The Lusternik-Schnirelmann category   (LS-category in short) of a smooth  manifold $M$ is an invariant putting in relation the topological complexity of $M$ with the behavior of smooth functions defined on $M$. In particular, when the manifold is compact, it is a lower bound of the number of critical points for any smooth function. This numerical invariant can also be defined for a topological space $X$ and reveals itself as a homotopy type invariant. 
The LS-category $\cat(X)$ is the least natural number $n$ such that there exists a cover of $X$ by $n+1$ open sets, each of them being contractible to a point inside $X$.

\smallskip
H. Colman and the second author (\cite{CMV01} and \cite{CMV02}, see also \cite{Col98}) have adapted this definition  to the case of a foliated manifold $(M,\cF)$. Indeed, they introduced two invariants, one which refers to the transverse structure and a second one, the \emph{tangential category}, $\cat_\cF(M)$,  see \defref{def:foliatedLS} for a precise statement. Our paper is concerned with this tangential category, also investigated by W.~Singhof and E.~Vogt \cite{SV03}. These two authors prove that the tangential category is less than or equal to dim $\cF + 1$, a result similar to the fact that the dimension of a CW-complex plus 1 is an upper bound for its LS-category. 
Some other results of \cite{SV03} will be crucial for our work and we will quote them precisely in the last section. 
In \cite{C-H02}, H. Colman and S. Hurder  prove, for instance, that the nilpotency index of the reduced filtered cohomology is a lower bound of the tangential category.

\smallskip
In \cite{Hur03}, S. Hurder says that the homotopy-theoretic interpretation of $\cat_\cF(M)$ corresponding to the Whitehead and Ganea definitions of category is ``one of the most important open problems in the subject''.
In the present work we give for the tangential category  several equivalent definitions inspired by the Whitehead procedure (using the fat wedge), and by the Ganea construction (using Milnor's  classifying spaces). As in the classical case, we do that in the topological 
framework, considering the most general setting of topological spaces, endowed with a partition. We call them \emph{stratified spaces} and denote by $\Tops$ the category of  stratified spaces and
stratified maps.


\smallskip
Before giving the content of the paper, we present {our general
 strategy}. Recall that the LS-category of a space $X$ uses open sets that
are contractible to a point inside $X$. If $X$ is path-connected, one can use
contractions on a specific point and work with pointed tools as fat wedges and
Ganea fibrations. These types of constructions are also present in our work
and, at a first look, we should distinguish some transverse structure
for a foliation. But we cannot do this, because the definition of tangential
LS-category is made of tangentially contractible open sets without specific
transverse set. We solve this contradiction by taking an adequate definition
of topological transverse set and the infimum on all transverse sets. It is a
priori unclear that this process gives the tangential LS-category of a
foliation and \secref{sec:LSfoliation} contains a proof of that fact.

\smallskip
The second main point is the equality between the invariants  coming from the Whitehead and Ganea constructions. As the first author showed in \cite{Doe93}, one has only to prove that $\Tops$ is a closed model category satisfying the Cube Lemma. Instead of that, we take here a shortcut:  the heart of the Cube Lemma needs only a structure of fibration category as it appears in the proof of \thmref{thm:ganeawhitehead}. Therefore, in \thmref{thm:stratemodel} of \secref{sec:closedmodelstratified}, we  prove  that $\Tops$ is a category of fibrations, which is sufficient for our purpose.


\smallskip
In \secref{sec:whiteheadganea}, we  consider a stratified pair $(X,A,\cF_X)$, where $(X,\cF_X)$ is a stratified space and $A$ a subset of $X$. 
We introduce the Whitehead and Ganea constructions which give the notions of Whitehead and Ganea category, respectively denoted 
$\Whcat(X,A,\cF_X)$ and $\Gcat(X,A,\cF_X)$. We prove their equality in \thmref{thm:ganeawhitehead}.
These constructions are adaptations of those by M. Clapp and D. Puppe \cite{Cl-Pu86} in the case of topological spaces.
In a third step, we want a notion of LS-category for $(X,A,\cF_X)$ using open sets. An open subset $U$ of $X$ is said \emph{$A$-categorical} if there is a stratified homotopy defined on $(U,\cF_U)$ between the identity on $U$ and a map with values in $A$, see \defref{def:Acategorique}. (The  induced stratification $\cF_U$ is composed of the connected components of the intersections $U\cap S$, $S\in \cF_X$.) This notion of $A$-categorical open set brings a definition of open LS-category, $\Ocat(X,A,\cF_X)$, in the usual way. When $A$ is a $B$ stratified neighborhood deformation (in short $B$-SND, see \defref{def:snr}) and $X$ a normal space, we prove that
$\Ocat(X,B,\cF_X)\leq\Whcat(X,A,\cF_X)\leq\Ocat(X,A,\cF_X)$, see \thmref{thm:equivalence3}.

\smallskip
In \secref{sec:LSstratified}, we introduce the \emph{transverse subsets} which are the key for the comparison with the tangential LS-category of foliations. If $(X,\cF_X)\in\Tops$, a subset $A$ of $X$ is called transverse to $\cF_X$ if the intersection $A\cap S$ is at most countable for all strata $S\in\cF_X$. The transverse subsets to a foliation own this property, see \cite{SV03}. In the definition of the tangential LS-category of foliations, there is no predetermined transverse set $A$, so we define $\Ocat(X,\cF_X)$ as the infimum of the integers $\Ocat(X,A,\cF_X)$ when $A$ is transverse to $\cF_X$
(analogously for $\Whcat(X,A,\cF_X)$ and $\Gcat(X,A,\cF_X)$).
In  \thmref{thm:invarianceOcat}, we prove that $\Ocat(X,\cF_X)$ is a homotopy invariant in $\Tops$. 

\smallskip
Finally, in \secref{sec:LSfoliation}, we consider a smooth closed manifold with a $C^1$-foliation, $(M,\cF_M)$, and prove the equalities:
$$ \cat_\cF(M)=\Ocat(M,\cF_M)=\Gcat(M,\cF_M)=\Whcat(M,\cF_M),$$
see \thmref{thm:tangentialopen}. We end with the example  of the Reeb foliation on the torus $T^2$ by giving a transverse set $A$ that is an $A$-SND.

\medskip
The second author would like to thank the University of Sciences and Technology of Lille and the first and third authors the University of Santiago de Compostela for their hospitality during the realization of this work.

\section{The category of stratified spaces}\label{sec:closedmodelstratified}
A \emph{stratified space} $(X,\cF_X)$ is a topological
space together with a partition $\cF_X$ whose
elements
$S\subseteq X$ are path-connected subspaces. We denote by
$\sim_{\cF_X}$ (or $\sim_X$ if there is no ambiguity) the equivalence relation associated to the
partition $\cF_X$ of $X$ and call \emph{strata} the elements $S$ of  $\cF_X$. A stratified map between stratified spaces
is a continuous map compatible with the equivalence
relations.

\smallskip
Let \Top~be the category of topological spaces and continuous maps 
and $\Tops$ be the category of stratified spaces and stratified maps.
We first observe that the category \Tops~has finite direct and inverse limits. 
This is obvious from the corresponding limits in \Top.

\smallskip
The stratification on a pull-back $(X,\cF_X)\times_{(B,\cF_B)} (Y,\cF_Y) = (P,\cF_P)$ is defined by $(x,y) \sim_P (x',y')$  if, and only if,  \\
\centerline{$x\sim_X x'$ \text{ and } $y\sim_Y y'$.}

\smallskip
If we have a commutative diagram in \Tops:
$$\xymatrix{(Z,\cF_Z)\ar[d]_i\ar[rr]^j&&(Y,\cF_Y)\ar[d]^g\\
(X,\cF_X)\ar[rr]^f&&(B,\cF_B)}$$
then the map induced by the pull-back in $\Top$, \\
\centerline{$(i,j): Z \to P = X\times_B Y$,}
is a stratified map in $\Tops$, because if $z \sim_Z z'$, then $i(z) \sim_X i(z')$ and $j(z) \sim_Y j(z')$, so $(i,j)(z) = (i(z),j(z)) \sim_P (i(z'),j(z')) = (i,j)(z')$.

\smallskip
The stratification on a push-out $(X,\cF_X)\vee_{(B,\cF_B)} (Y,\cF_Y) =
(S,\cF_S)$ 
of $f: (B,\cF_B) \to (X,\cF_X)$ and $g: (B,\cF_B) \to (Y,\cF_Y)$
is defined by
\begin{itemize}
\item $x \sim_S x'$ if, and only if, 
\begin{itemize}
\item \emph{either}: $x\sim_X x'$, 
\item \emph{or}:  $x  \sim_X f(b)$  and  $x' \sim_X f(b')$  with  $g(b) \sim_Y g(b')$, (\dag) 
\end{itemize}
\item a similar formula for $y \sim_S y'$,
\item and $x \sim_S y$ if, and only if, \\
\centerline{$x  \sim_X f(b)$  and  $y \sim_Y g(b)$ (\ddag).}
\end{itemize}

\smallskip
If we have a commutative diagram in \Tops:
$$\xymatrix{(B,\cF_B)\ar[d]_f\ar[rr]^g&&(Y,\cF_Y)\ar[d]^j\\
(X,\cF_X)\ar[rr]^i&&(Z,\cF_Z)}$$
then the map induced by the push-out in $\Top$, \\
\centerline{$[i,j]: S = X\vee_B Y \to Z$,}
is a stratified map in $\Tops$, because if $x \sim_S x'$, with condition (\dag), then
$i(x) \sim if(b) = jg(b) \sim jg(b') = if(b') \sim i(x')$ and if $x \sim_S y$, with condition (\ddag), then $i(x) \sim if(b) = jg(b) \sim j(y)$.  


\smallskip
Denote by $I$ the interval $[0,1]$ endowed with the trivial partition $\cF_I=\{[0,1]\}$ and define a notion of homotopy in $\Tops$ as follows.

\begin{definition}\label{def:homotopystrate}
Two stratified maps, $f,g\colon (X,\cF_X)\to (Y,\cF_Y)$, are called \emph{$\cS$-homotopic} if there exists a stratified map
$H\colon (X,\cF_X)\times (I,\cF_I)\to (Y,\cF_Y)$
such that $f=F(-,0)$ and $g=F(-,1)$. We denote this relation by $f\simeq_\cS g$.
\end{definition}

Recall that \Top~can be endowed with two structures of closed model category. The first one (see \cite{Qui67}) has for fibrations the Serre fibrations and for weak equivalences the classical weak equivalences. The second one (see \cite{Str72}) has for  fibrations
the Hurewicz fibrations and for weak equivalences the homotopy
equivalences. For the study of \Tops, we consider a situation ``\`a la Hurewicz" and prove that \Tops~is a fibration category, in the sense of Baues
(see \cite[I {\S}1a]{Bau89}), which is sufficient for our objective.

\begin{theorem}\label{thm:stratemodel}
The category \Tops~is a fibration category with

-- \emph{fibrations,} the stratified maps $f\colon (X,\cF_X)\to (Y,\cF_Y)$ having the $\cS$-homotopy lifting property;

-- \emph{weak equivalences,} the stratified maps $f\colon (X,\cF_X)\to (Y,\cF_Y)$ such that there exists a stratified map $g\colon (Y,\cF_Y)\to (X,\cF_X)$ with $f\circ g\simeq_\cS \id_Y$ and $g\circ f\simeq_\cS \id_X$.
\end{theorem}
Note that the full subcategory of stratified spaces $(X,\cF_X)$~with $\cF_X$ formed of the path-connected components of $X$ equals \Top~with its structure of fibration category.
In particular, a weak equivalence in \Tops~is a weak equivalence in \Top.
\begin{proof} As the proof is an adaptation of \cite{Str72} to the stratified case, we only sketch it.
Recall that a fibration category is a category satisfying the axioms (F1), (F2), (F3) and (F4) of \cite{Maj00}.

\smallskip
(F1) The isomorphisms are weak equivalences and also fibrations. Given two maps \\
\centerline{$A \xrightarrow{f} B \xrightarrow{g} C$}
if any two of $f, g, gf$ are weak equivalences, then so is the third. The composite of fibrations is a fibration. 
The verification of these properties is immediate.

\smallskip
(F2) Consider a pull-back in \Tops:
$$\xymatrix{(P,\cF_P)\ar[d]_q\ar[rr]^g&&(E,\cF_E)\ar[d]^p\\
(Y,\cF_Y)\ar[rr]^f&&(B,\cF_B)}$$

We have to show that if $p$ is a fibration of \Tops, then so is $q$. 
Indeed, from the universal property, the map $q$ has the $\cS$-homotopy lifting property and therefore is a fibration of \Tops.

We have also to show that if $p$ is a trivial fibration of \Tops, then so is $q$. 
As $p$ is a weak equivalence, there exists $\sigma: (B,\cF_B)\to (E,\cF_E)$ such that $p\circ \sigma\simeq_\cS \id_B$ and $\sigma\circ p\simeq_\cS\id_E$. As $p$ is  a fibration, we may suppose that $p\circ \sigma=\id_B$. From the universal property, we deduce the existence of $\mu\colon (Y,\cF_Y)\to (P,\cF_P)$ such that $g\circ \mu=\sigma\circ f$ and $q\circ \mu=\id_Y$. We have to show that $\mu\circ q\simeq_\cS \id_P$.  For that, we follow the proof in the case of the category \Top~\cite{Str72}, adding the ad~hoc argument for the stratification. We briefly recall the main steps for the convenience of the reader:
\begin{itemize}
 \item First, one knows (see \cite{Dol63}) that the homotopy $H$ between $\sigma\circ p$ and the identity on $E$ can be chosen such that $p\circ H(-,t)=p$ for any $t\in I$. 
 \item With that homotopy $H$, we check the commutativity of the following diagram, where $\proj$ means a projection onto one factor.
$$\xymatrix{
&P\times I\ar@{..>}[ld]_G\ar[rr]^{g\times \id}
\ar[dd]|(.49)\hole_(.6){q\times \id}
&&E\times I\ar[dd]^{p\times \id}\ar[ld]_H\ar[dd]\\
P\ar[rr]^(.7)g\ar[dd]_{q}&&
E\ar[dd]_(.66){p}&\\
&Y\times I\ar[rr]|(.49)\hole^(.66){f\times \id}
\ar[ld]_{\proj}&&B\times I\ar[ld]^{\proj}\\
Y\ar[rr]_f&&B&
}$$  
From the universal property, there exists 
$G\colon (P,\cF_P)\times (I,\cF_I)\to (P,\cF_P)$ such that 
$g\circ G= H\circ (g\times \id)$ and $q\circ G=\proj\circ (q\times \id)$. This map $G$ satisfies $G(-,1)=\id$ and $G(-,0)=\mu\circ q$ as wanted. 
\end{itemize}

On the other hand, it is easily checked that any object is fibrant, which implies, as quoted by H.~Baues (see \cite[Lemma 1.4 page 7]{Bau89}), that the category \Tops~is `proper', i.e. if we assume that $f$ is a weak equivalence, so is $g$.

\smallskip
(F3) Let $f\colon (X,\cF_X)\to (Y,\cF_Y)$   be a map in \Tops. We have to write this map as $f = p\circ s$ where $s$ is a weak equivalence in \Tops~and $p$ is a fibration in \Tops. Recall that $Y^I$ denotes the space of continuous maps $\omega\colon I\to Y$. 
Here we denote by $(Y,\cF_Y)^I$ the set of continuous maps, $\omega\colon I\to Y$, such that $\omega(I)$ is included in an element $S$ of $\cF_Y$, together with the partition formed of the subsets $S^I$ with $S\in\cF$.
The proof of the axiom (F3) follows the classical way \cite{Str72}, substituing $Y^I$ by $(Y,\cF_Y)^I$. The different steps of the proof are:
\begin{itemize}
 \item The map $p_1\colon (Y,\cF_Y)^I\to (Y,\cF_Y)$ defined by $p_1(\omega)=\omega(1)$ is a fibration in \Tops.
 \item The map $s_1\colon (Y,\cF_Y)\to (Y,\cF_Y)^I$, sending $ y$ to the constant path $\hat y$, is a $\cS$-homotopy inverse of $p_1$.
 \item If we denote by $P_f$ the pull-back of $p_1$ and $f$, the map $p_0\colon P_f\to Y$ defined by $p_0(\omega,x)=\omega(0)$ is a fibration. 
 \item The map $s_f\colon (X,\cF_X)\to P_f$, given by $ x \mapsto (\widehat{f(x)},x)$, is a weak equivalence of \Tops~that gives the desired decomposition of $f$ as $f=p_0\circ s_f$.
\end{itemize}

\smallskip
(F4)  We check easily that every trivial fibration of \Tops~has a section, which implies (F4).
\end{proof}

\section{Lusternik-Schnirelmann category of stratified pairs} \label{sec:whiteheadganea}

\begin{definition}\label{def:canonicaltops}
Let $(X,\cF_X)\in \Tops$ be a stratified space and $A$ be a topological subspace of $X$. The \emph{induced stratification} $(A,\cF_A)$   has for strata the connected components of the intersections $A\cap S$, with $S\in \cF_X$. The pair $(A,\cF_A)$ is called \emph{a stratified subspace} of $(X,\cF_X)$ and $(X,A,\cF_X)$ a \emph{stratified pair.}
\end{definition}

If $f\colon (X,\cF_X)\to (Y,\cF_Y)$ is a stratified map, we observe from the previous definition, that the restriction 
$f\colon (X,\cF_X)\to (f(X),\cF_{f(X)})$ is also stratified.

\subsection{Ganea construction: the invariant $\Gcat(X,A,\cF)$}

\begin{definition}\label{def:ganea}
The $i$th  \emph{Ganea space} $G_i(X,A)$ of a stratified pair $(X,A,\cF_X)$ is the set of $(i+1)$-uples $(\alpha_1,\ldots,\alpha_{i+1})$ of paths in $(X,\cF_X)^I$, such that
$\alpha_1(1)=\dots=\alpha_{i+1}(1)$ and there exists $k$, $1\leq k\leq i+1$, with
$\alpha_k(0)\in A$.
\end{definition}

We observe from this definition that, for any element $(\alpha_1,\ldots,\alpha_{i+1})$ of $G_i(X,A)$, there is an associated element $S$ of the  partition $\cF_X$ such that $\alpha_j(I)\subset S$, for all $j=1,\ldots,i+1$. 

\smallskip
We put on $G_i(X,A)$ the stratification induced by the product stratification and the stratification of $X^I$, and denote by $G_i(X,A,\cF_X)$ this stratified space. The $i$th \emph{Ganea fibration} $g_i\colon G_i(X,A,\cF_X)\to (X,\cF_X)$ is defined by
$g_i(\alpha_1,\ldots,\alpha_{i+1})=\alpha_1(1)$.

\begin{definition}\label{def:LSganea}
The \emph{Ganea LS-category} of a stratified pair $(X,A,\cF_X)$
is the least integer $n$ such that the $n$th {Ganea fibration}
$g_n\colon G_n(X,A,\cF_X)\to (X,\cF_X)$
has a stratified section. We denote it by  $\Gcat(X,A,\cF_X)$.
\end{definition}

\subsection{Whitehead construction: the invariant $\Whcat(X,A,\cF_X)$}
\begin{definition}\label{def:whitehead}
The $i$th \emph{fat wedge} of a stratified pair $(X,A,\cF_X)$
is the subset $T_i(X,A)$ of the product $X^{i+1}$ formed of the elements 
$(x_1,\ldots,x_{i+1})$ such that there exists $k$, $1\leq k\leq i+1$, with $x_k\in A$.

We endow it with the stratification induced from the stratification of the product $(X,\cF_X)^{i+1}$ and denote this stratified space by $T_i(X,A,\cF_X)$.
\end{definition}

\begin{definition}\label{def:LSwhitehead}
The \emph{Whitehead LS-category} of a stratified pair
$(X,A,\cF_X)$ is the least integer $n$ such that the diagonal 
$\Delta: (X,\cF_X) \to (X,\cF_X)^\npu$ factors up to $\cS$-homotopy through the inclusion map $t_n: T_n(X,A,\cF_X) \to (X,\cF_X)^\npu$.
We denote it by $\Whcat(X,A,\cF_X)$.
\end{definition}

If $A=\left\{\ast\right\}$, they coincide with the classical ones, see \cite{CLOT03}, \cite{Doe93}, \cite{Cl-Pu86} or \cite{FVA03} for more details.

\subsection{Open LS-category: the  invariant $\Ocat(X,A,\cF_X)$}\label{subsec:opencat}

In this paragraph, we provide a third definition using open sets.
First we precise the notion of categorical subsets in the context of stratified spaces.

\begin{definition}\label{def:Acategorique}
Let  $(X,A,\cF_X)$ be a stratified pair.  A non-empy subset $U$ of $X$ is said to be \emph{$A$-categorical} if there is a stratified homotopy
$H\colon (U,\cF_U)\times (I,\cF_I)\to (X,\cF_X)$ such that $H(x,0)=x$ and $H(x,1)\in A$ for any $x\in U$. We call $H$ \emph{a stratified deformation of $U$ into $A$.}
\end{definition}

Observe that, as in the classical case, the sets $H(U,t)$, with $t\in I$, are not necessarily contained in $U$. 
With the previous definition, we introduce the \emph{open LS-category} of the stratified pair $(X,A,\cF_X)$.

\begin{definition}\label{def:opencat}
The \emph{open LS-category} of a stratified pair
$(X,A,\cF_X)$  is the least integer $n$ such that there exists a covering of $X$ by $(n+1)$ open sets which are $A$-categorical.  We denote it by  $\Ocat(X,A,\cF_X)$.  
\end{definition}

We write also $\Ocat(X,A,\cF_X)=\infty$ if such a covering does not exist.
Observe   that $\Ocat(X,X,\cF_X)=0$ and that the finiteness of $\Ocat(X,A,\cF_X)$ implies that $A$ cuts each stratum.

\begin{example}\label{exam:O(2)}
Examples with an infinite value for \Ocat can be obtained easily. For space $X$, we take the plane $\R^2$ stratified by the orbits of the action of the rotation group $\SO(2)$. Let $A$ be the half-ray $Ox$. We observe that any open subset $U$ containing the singular orbit $O$ contains also a circle $C$ which cannot be contracted to some point of $A\cap C$ by a stratified deformation of $U$. Therefore $\Ocat(\R^2,Ox,\cF_{\R^2})=\infty$.
\end{example}

The next easy result will be used in the proof of \thmref{thm:invarianceOcat}.

\begin{proposition}\label{prop:ocatf(A)}
Let  $(X,A,\cF_X)$ be a stratified pair.
If $f\colon (X,\cF_X)\to (Y,\cF_Y)$ is a stratified map with a right $\cS$-homotopical inverse $g$, then\\
\centerline{$\Ocat(Y,f(A),\cF_Y)\leq \Ocat(X,A,\cF_X)$.}
\end{proposition}

\begin{proof}
Let $U\subseteq X$ be an $A$-categorical open set and $H\colon (U,\cF_U)\times (I,\cF_I)\to (X,\cF_X)$ be a stratified deformation of $U$ into $A$. On the open set  $V=g^{-1}(U)$, we define a stratified deformation
$f\circ H\circ (g\times {\rm id}_I)\colon (V,\cF_V)\times (I,\cF_I)\to (Y,\cF_Y)$ into $f(A)$. The inequality follows.
\end{proof}

\subsection{Comparison of the three invariants}

\begin{theorem}\label{thm:ganeawhitehead}
For any stratified pair  $(X,A,\cF_X)$, we have the equality\\
\centerline{$\Whcat(X,A,\cF_X)= \Gcat(X,A,\cF_X)$.}
\end{theorem}

\begin{proof}
We first define a stratified map $\varepsilon_n\colon G_n(X,A,\cF_X)\to T_n(X,A,\cF_X)$ by $\varepsilon_n(\alpha_1,\ldots,\alpha_{n+1})=
(\alpha_1(0),\ldots,\alpha_{n+1}(0))$. The following diagram is clearly homotopy commutative:
$$\xymatrix{G_n(X,A,\cF_X)\ar[rr]^{g_n}\ar[d]_{\varepsilon_n}
&&(X,\cF_X)\ar[d]^\Delta\\
T_n(X,A,\cF_X)\ar[rr]_{t_n}&&(X^{n+1},\cF_{X^{n+1}})
}$$
The proof is reduced to the fact that this diagram  is a homotopy pull-back. For that, we need to determine the associated fibration to the diagonal $\Delta$,  which is obtained from the pull-back
$$\xymatrix{
(P,\cF_P)\ar[r]\ar[d]&
(X^{n+1},\cF_{X^{n+1}})^I\ar[d]\\
(X,\cF_X)\ar[r]_-\Delta&(X^{n+1},\cF_{X^{n+1}})
}$$
By construction, $(P,\cF_P)$ is the space of
$(n+1)$-uples $(\alpha_1,\ldots,\alpha_{n+1})$ of paths in $(X,\cF_X)^I$, such that
$\alpha_1(1)=\dots=\alpha_{n+1}(1)$. The evaluation map to $0$ of these paths in $X^{n+1}$ is a fibration equivalent to the diagonal. So if we take the pull-back of this map and of the map
$t_n: T_n(X,A,\cF_X)\to(X^{n+1},\cF_{X^{n+1}})$, we get the homotopy pull-back of $t_n$ and of the diagonal. But this 
is exactly the Ganea space $G_n(X,A,\cF_X)$.
\end{proof}

We prove now the equivalence between the Whitehead and the open-set definitions, under some hypotheses, as in the classical case where one needs the existence of a contractible neighborhood of the base point. The corresponding notion is given by the following definition.

\begin{definition}\label{def:snr}
Let $(X,\cF_X)$ be a stratified space and let $A, B$ be subsets of $X$.
We say that $A$ is a \emph{$B$ stratified neighborhood deformation} (in short 
$B$-SND) if $A$ has some open neighborhood which is a $B$-categorical set.
\end{definition}

\begin{theorem}\label{thm:equivalence3}
Let $(X,\cF_X)$ be a stratified space and $A$ be a subspace of $X$.
\begin{enumerate}
 \item If the space $X$ is normal, then we have the inequality:
\item[] \centerline{$\Whcat(X,A,\cF_X)\leq \Ocat(X,A,\cF_X)$.}
 \item If $A$ is a $B$-SND, then we have the inequality:
\item[] \centerline{$\Ocat(X,B,\cF_X)\leq \Whcat(X,A,\cF_X)$.}
 \end{enumerate}
\end{theorem}

\begin{proof}
1)  Suppose $\Ocat(X,A,\cF_X)\leq n$. There is a covering of $X$ by open sets, $U_0,\ldots, U_n$,  and stratified deformations
$H_i\colon (U_i,\cF_{U_i})\times (I,\cF_I)\to (X,\cF_X)$ into $A$, for $i=0,\ldots,n$. As $X$ is normal, there exists a covering of $X$ by open sets, $W_0,\ldots,W_n$, such that $\overline{W}_i\subset U_i$, for $i=0,\ldots,n$. For any $i$, we choose an Urysohn function $\varphi_i\colon X\to I$ such that $\varphi_i(x)=1$ if $x\in \overline{W}_i$ and $\varphi_i(x)=0$ if $x\notin U_i$.
We define now a continuous stratified map $\hat{H}_i\colon (X,\cF_X)\times (I,\cF_I)\to (X,\cF_X)$ by:
$$\hat{H}_i(x,t)=\left\{
\begin{array}{cl}
     H_i(x,\varphi_i(x)t)  & \mathrm{if\  } x\in U_i, \\
     x  & \mathrm{otherwise}.\\
\end{array}
\right.$$
We collect these maps in a continuous stratified map
$H\colon (X,\cF_X)\times (I,\cF_I) \to (X,\cF_X)^{n+1}$ defined by
$H(x,t)=(\hat{H}_0(x,t),\ldots,\hat{H}_n(x,t))$. Observe that we have $H(x,0)=(x,\ldots,x)=\Delta(x)$.

Set $r(x)=H(x,1)$. Since the $W_i$'s are a covering of $X$, for any point $x\in X$, there is a $W_j$ with $x\in W_j$. By definition of $\hat{H}_j$, we have $\hat{H}_j(x,1)=H_j(x,1)\in A$. From the construction of the fat wedge, we have $r(X)\subset T_{n}(X,A,\cF_X)$ and $r$ is a lifting up to $\cS$-homotopy (by the $\cS$-homotopy $H$) of the diagonal. By definition, we get
$\Whcat(X,A,\cF_X)\leq n$.

\smallskip
2) Suppose $\Whcat(X,A,\cF_X)\leq n$. By definition, there are a stratified map
$r\colon (X,\cF_X)\to T_n(X,A,\cF_X)$ and a stratified homotopy $H\colon (X,\cF_X)\times (I,\cF_I)\to (X,\cF_X)^{n+1}$ between the diagonal $\Delta$ and the composite $t_n\circ r$, see \defref{def:LSwhitehead}. As $A$ is a $B$-SND, there exists also an open set $N$, $A\subset N$, and a stratified homotopy $G\colon (N,\cF_N)\times (I,\cF_I)\to (X,\cF_X)$ of $N$ with $G(x,1)\in B$.

Let $p_i\colon X^{n+1}\to X$ be the $(i+1)$th projection, $0\leq i \leq n$, we set $h_i= p_i\circ t_n\circ r$ and
$U_i=h_i^{-1}(N)$. Then, since
$r(X)\subset T_{n}(X,A,\cF_X) = \bigcup_{i=0}^{n}p_i^{-1}(A)$ we have $X=\bigcup_{i=0}^{n}U_i$. The $U_i$'s are a covering of $X$. They are $B$-categorical with the homotopy $H_i\colon U_i\times I\to X$ defined by:
$$H_i(u,t)=\left\{
\begin{array}{cl}
     p_iH(u,2t)  & \mathrm{if\  } 0\leq t \leq 1/2, \\
    G(h_i(u),2t-1)  & \mathrm{if\ }  1/2\leq t \leq 1.\\
\end{array}
\right.$$
It is routine to check that $H_i$ is a continuous homotopy between the identity and a map with values in $B$. As $H$, $p_i$, $G$ and $h_i$ are stratified maps, $H_i$ is stratified also and we have proved
$\Ocat(X,B,\cF_X)\leq n$.
\end{proof}

Let $X$ be path-connected with the trivial stratification and let $U$ be a
categorical open set in $X$. For any point $u\in U$ and any $x\in X$, the
subset $\{u\}$ is an $\{x\}$-SND. 
But, in the general case of a stratified space, being $B$-SND depends on the stratification but also depends on $B$. This explains why we need the notion of transverse sets, introduced in the next section.

\section{Lusternik-Schnirelmann  category of stratified spaces}\label{sec:LSstratified}

First, we define the \emph{transverse subsets} of a stratified space $(X,\cF_X)$ which play a fundamental role for the comparison of our previous invariants  with the tangential LS-category of foliations.

\begin{definition}\label{def:transverse}
Let $(X,\cF_X)$ be a stratified space. A subspace $A\subseteq X$ is \emph{transverse to the stratification} if, for any stratum $S\in\cF_X$, the set $A\cap S$ is at most countable. We denote this property by $A \transverse \cF_X$.
\end{definition}

As it is proved by W. Singhof and E. Vogt, \cite{SV03}, the transverse subspaces of a foliation are transverse in the sense of Definition \ref{def:transverse}. We observe now that, in the \defref{def:foliatedLS} of the tangential category of a foliation, the transverse space which receives the stratified deformation is not predetermined. We have therefore to take in account all the transverse subspaces associated to tangential deformations.

\begin{definition}\label{def:OcatX}
Let  $(X,\cF_X)$ be a stratified space. The \emph{open LS-category} of
$(X,\cF_X)$  is the infimum of the integers $\Ocat(X,A,\cF_X)$, when $A$ runs along the transverse subsets to $\cF_X$:
$$\Ocat(X,\cF_X)={\rm Inf\;}\left\{\Ocat(X,A,\cF_X)\mid A\transverse \cF_X\right\}.$$
Similarly, we define the \emph{Whitehead LS-category and the Ganea LS-category} of the stratified space $(X,\cF_X)$ by taking the infimum along the transverse subsets. We denote them by
$\Gcat(X,\cF_X)$ and  $\Whcat(X,\cF_X)$ respectively.
\end{definition}

The lower bound gives a homotopy invariant. 

\begin{theorem}\label{thm:invarianceOcat}
In the fibration category $\Tops$, the integer $\Ocat(-,-)$ is a homotopy invariant.
\end{theorem}

\begin{proof}
Let $f\colon (X,\cF_X)\to (X',\cF_{X'})$ be a stratified map with homotopical inverse
$g$ in \Tops. Let~$\sim$ (resp.~$\sim'$) be the equivalence relation generated by $\cF_X$ (resp. $\cF_{X'}$). The stratum containing $x\in X$ (resp. $f(x)\in X'$) is denoted by $S_x$ (resp. $S'_{f(x)}$).
The proof is divided in several steps.

\smallskip
1) First, one proves that \emph{$f$ induces a homeomorphism $\bar{f}\colon X{/\!\!\sim}\to X'{/\!\!\sim'}$ between the quotient spaces.} 
If $y\in S'_{f(x)}$, there is a path $\gamma\colon [0,1]\to S'_{f(x)}$ joining $y$ and $f(x)$. Therefore $g(y)$ and $g(f(x))$ are connected by a path included in one stratum. As $g\circ f\simeq_\cS {\rm id}_X$, the points $g(f(x))$ and $x$ are in the same stratum. This implies $g(S'_{f(x)})\subseteq S_x$ and
$\bar{g}\circ \bar{f}={\rm id}_{X/\sim}$. 
With a similar argument, one has $\bar{f}\circ \bar{g}={\rm id}_{X'/\sim'}$.

\smallskip
2) Suppose $f(S)\subseteq S'$ where $S$ and $S'$ are strata of $\cF_X$ and
$\cF_{X'}$ respectively. Let $A\subset X$. If $f(a)\in f(A)\cap S'$, then $g(f(a))\in S$, from Part~1, and $a\in S$, because $a$ and $g(f(a))$ are in the same stratum. This implies $f(A)\cap S'\subseteq f(A\cap S)$. As the reverse inclusion is obvious, one has
$f(A\cap S)=f(A)\cap S'$.

\smallskip
3) Let $A\subset X$ such that $f(A)$ is not transverse to $\cF_{X'}$, i.e. the set $f(A)\cap S'$ is uncountable for some stratum $S'\in\cF_{X'}$. Let $S\in \cF_X$ with $g(S')\subseteq S$. From the equality $f(A\cap S)=f(A)\cap S'$, we deduce that $A\cap S$ is uncountable. We have proved that $A\transverse \cF_X$ implies $f(A)\transverse \cF_{X'}$.

\smallskip
4) \propref{prop:ocatf(A)} implies the inequality
$${\rm Inf\;}\left\{\Ocat(X,A,\cF_X)\mid A\transverse \cF_X\right\}\geq
{\rm Inf\;}\left\{\Ocat(X',A',\cF_{X'})\mid A'\transverse \cF_{X'}\right\}.$$
As the opposite inequality follows by symmetry, we have proved the equality 
$\Ocat(X,\cF_X)=\Ocat(X',\cF_{X'})$.
\end{proof}

\begin{remark}\label{rem:totallydisconnected}
In the previous proof, we have used, in a fundamental way, the fact that the image,
by a weak equivalence of \Tops,
of a transverse set  of $(X,\cF_X)$ is a transverse set of $(X',\cF_{X'})$. With the next example, we observe that this is not true for the subsets $A$ of $(X,\cF_X)$ such that $A \cap S$ is totally disconnected for any $S\in \cF_X$.
For similar reasons, the property $A \cap S$ of topological dimension~0 does not fit also, see \cite[Theorem Page~302]{Kur}.
\end{remark}
\begin{example}\label{exam:cylinder}
Let $f\colon X=S^1\times \R\to X'=S^1$ be the  projection on the first factor. We put on $X$ and $X'$ the trivial stratification with only one stratum. 
Let $\varphi\colon S^1\to \R$ be an application that is discontinuous at each point and let $A$ be the graph of $\varphi$.  We see easily that $f$ is a weak equivalence of \Tops, that $A$ is totally disconnected and $f(A)=X'$ is not.
\end{example}

\section{Tangential LS-category of foliations}\label{sec:LSfoliation}

The {\it tangential LS-category} of a foliated manifold was introduced by H. Colman and the second author in \cite{CMV02}, see also \cite{Col98}. In this section we compare it with our previous definitions of  stratified LS-categories; all our manifolds are $C^\infty$-manifolds.

\smallskip
Let $M$ be a foliated manifold of class $C^0$. The leaves of the foliation form a partition $\cF_M$ of $M$. Moreover the definitions of foliated continuous maps and tangential continuous homotopies correspond exactly to our definitions of stratified maps and stratified homotopies. The induced foliation on an open set of $M$ coincides also with our induced stratification of
\defref{def:canonicaltops}.
In the case of a foliation of class $C^\infty$ on a closed manifold, recall from  \cite[Theorem 4.1]{SV03} that the use of continuous maps or of $C^\infty$-maps does not modify the tangential LS-category. 

\subsection{Definition of tangential LS-category}

\begin{definition}\label{def:foliatedcategorical}
An open set $U\subset M$ (endowed with the induced foliation $\cF_U$) is \emph{tangentially categorical} if there exists a tangential homotopy $H \colon U\times I \to M$ such that $H(-,0)$ is the inclusion and $H(-,1)$ is constant along the leaves of $U$, that is $H(-,1)$ maps each leaf of $\cF_U$ onto some point. We call $H$ a \emph{tangential contraction}.
\end{definition}

\begin{example}\label{exam:noextension}
In general  a tangential contraction defined on the categorical open set $U$ cannot be extended to the whole manifold $M$. For instance, let $M$ be the punctured plane $\R^2\backslash\{(0,0)\}$, foliated by the horizontal lines
$\R\timesÊ\{t\}$. Let $L$ be the closed leaf $L=\{(x,0)\colon x<0\}$ and take $U=M\backslash L$.
Then $U$ can be tangentially contracted to $A=\{1\}\times \R$, but the contraction cannot be extended to $M$.
\end{example}

\begin{definition} \label{def:foliatedLS} \cite{Col98, CMV02}
The {\it tangential LS-category} of the foliated manifold $(M,\cF_M)$ is the least integer $n$ such that there exists a covering of $M$   by $(n+1)$ open sets which are tangentially categorical. We denote it by $\catF(M)$.
\end{definition}

\subsection{Equivalence of the four invariants}
This paragraph is devoted to the proof of the equivalence between the tangential LS-category of a foliation and our definitions of the LS-category  of the corresponding stratified space. Let $(M,\cF_M)$ be a foliated manifold of class $C^0$ and dimension $p$.

\begin{lemma}\label{SV} Consider a tangentially categorical open set $U$ with its tangential contraction $H\colon U\times I \to M$. Then $H(U,1)$ is contained in a set $A(U)$ which is a finite or countable union of compact $(n-p)$-submanifolds with boundary, each one transverse to $\cF_M$. Moreover, for each leaf $L$ of $\cF_M$ the intersection $L\cap A(U)$ is a countable set.
\end{lemma}
\begin{proof}The first part is \cite[Lemma 1.1]{SV03}. The second part follows from the fact that each compact transverse submanifold can be covered by a finite number of adapted 
charts. Since its dimension is the codimension of the foliation, its intersection with a leaf  is a second 
countable manifold of dimension zero. 
\end{proof} 

\begin{theorem}\label{thm:tangentialopen} If $\cF_M$ is a $C^0$-foliation on the manifold $M$, we have the equality $\Ocat(M,\cF_M)=\catF(M)$. 
In the case of a $C^1$ foliation, the four invariants coincide:
$$\Gcat(M,\cF_M)=\Whcat(M,\cF_M)=\Ocat(M,\cF_M)=\cat_\cF(M).$$ 
\end{theorem}

\begin{proof} We first prove that $\cat_\cF(M) = \Ocat(M,\cF_M)$.

\smallskip
$\bullet$ Suppose $\cat_\cF(M)<\infty$. Let $U_0,\ldots,U_n$ be a covering of $M$ by tangentially categorical open sets. Denote by $H_i\colon U_i\times I\to M$ a tangential contraction of $U_i$ and by $A_i$ the set $A(U_i)$ containing  $H_i(U_i,1)$ as in Lemma \ref{SV}. Let $A=A_0\cup \cdots \cup A_n$. Then the set $A$ is  transverse to $\cF_M$ in the sense of  \defref{def:transverse}, and the open sets $U_i$ are $A$-categorical (\defref{def:Acategorique}). Hence $\Ocat(M,\cF_M)\leq n$.

\smallskip
$\bullet$ On the other hand, suppose $\Ocat(M,\cF_M)<\infty$. Since this infimum must be a minimum, there exists  a transverse subset $A$  such that $\Ocat(M,\cF_M) = \Ocat(M,A,\cF_M)$. Let $U$ be an $A$-categorical open set and
$H$ be the corresponding stratified homotopy. Then for each leaf $L\in \cF_M$, the intersection $L\cap A$ is a countable, thus totally disconnected set.
Hence the image by $H(-,1)$ of each connected component of $L\cap U$    is reduced to a point. Therefore, $U$ is tangentially categorical, which proves that $\cat_\cF(M)\leq n$.

Now, theorems \ref{thm:ganeawhitehead} and \ref{thm:equivalence3} (part 1) give:
$$\Gcat(M,\cF_M)=\Whcat(M,\cF_M)\leq\Ocat(M,\cF_M)=\cat_\cF(M).$$
So we just have to prove that $\Ocat(M,\cF_M)\leq\Whcat(M,\cF_M)$.

\smallskip
$\bullet$ If $\Whcat(M,\cF_M) = n$,  there exists $A$ such that 
$\Whcat(M,A,\cF_M) = n$ and $A\transverse \cF_M$. By Lemma 1.1 and Proposition 5.1 of \cite{SV03}, there is a subset $B$, such that $B\transverse \cF_M$ and $A$ is a $B$-SND. Now \thmref{thm:equivalence3} implies the inequality $\Ocat(M,\cF_M)\leq \Whcat(M,\cF_M)$.
\end{proof}

\begin{example}[Reeb foliation]\label{exam:reeb}
We ask if for an arbitrary foliation there exists a transverse subset $A$ such that we have the equality $\Ocat(M,\cF_M)=\Ocat(M,A,\cF_M)$ and that is an $A$-SND, i.e. a SND of itself.
Here we construct explicitly such $A$ for the example of the Reeb foliation on $T^2$, as described in \cite[Section 6.2]{CMV02}. 
Consider the following representation of the torus where the upper and lower arcs are identified. At the bottom, the left and right segments have also to be identified.
The strata are the horizontal lines.

\vskip-1.6cm
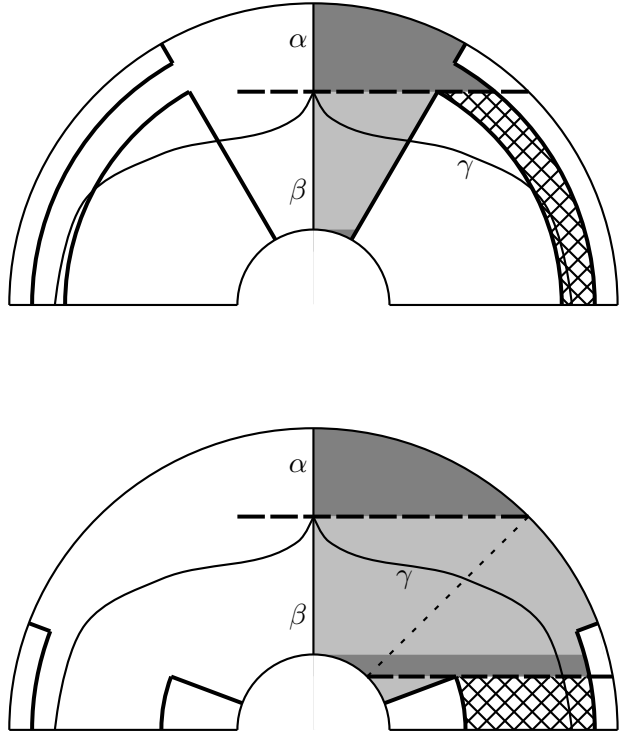
\centerline{\begin{pspicture}(-5,-1)(5,6)
\pscustom[fillstyle=solid, fillcolor=gray, linestyle=none]{%
  \psline(0,2.8284)(2.8284,2.8284)
  \psarc(0,0){4}{45}{90}
} %
\pspolygon[fillstyle=solid, fillcolor=lightgray, linestyle=none](0,0)(0,2.8284)(1.6330,2.8284)
\pspolygon[fillstyle=solid, fillcolor=gray, linestyle=none](0,1)(1,1)(1,0.7071)(0,0.7071)
\pscustom[fillstyle=solid, linestyle=none]{%
  \psline(0,0)(4,0)
  \psarc(0,0){4}{0}{60}
} %
\pscustom[fillstyle=solid, fillcolor=gray, linestyle=none]{%
  \psline(0,2.8284)(2.1056,2.8284)
  \psarc(0,0){3.7}{49.8572}{90}
} %
\pscustom[fillstyle=crosshatch, fillcolor=lightgray, linestyle=none]{%
  \psline(3.2660,0)(3.7,0)
  \psarc(0,0){3.7}{0}{49.8572}
  \psline(2.1056,2.8284)(1.6330,2.8284)
} %
\pscustom[fillstyle=solid, linestyle=none]{%
  \psline(0,0)(3.2660,0)
  \psarc(0,0){3.2660}{0}{60}
} %
\pscustom[fillstyle=solid, linestyle=none]{%
  \psarc(0,0){1}{0}{180}
  \psline(0,1)(0,0)
  \psline(0,0)(1,0)
} %
\psarc(0,0){4}{0}{180}
\psarc(0,0){1}{0}{180}
\psline(-4,0)(-1,0)
\psline(1,0)(4,0)
\psline(0,1)(0,4)
\pscurve(0,2.8284)(0.2,2.5)(2,2)(3,1.4)(3.4,0)
\pscurve(-0,2.8284)(-0.2,2.5)(-2,2)(-3,1.4)(-3.4,0)
\psline[linestyle=dashed, linewidth=1.5pt, dash=10pt 2pt](-1,2.8284)(2.8284,2.8284)
\psline[linewidth=1.5pt](0.5,0.8660)(1.6330,2.8284)
\psarc[linewidth=1.5pt](0,0){3.2660}{0}{60}
\psline[linewidth=1.5pt](1.85,3.2043)(2,3.4641) 
\psarc[linewidth=1.5pt](0,0){3.7}{0}{60}
\psline[linewidth=1.5pt](-0.5,0.8660)(-1.6330,2.8284)
\psarc[linewidth=1.5pt](0,0){3.2660}{120}{180}
\psline[linewidth=1.5pt](-1.85,3.2043)(-2,3.4641) 
\psarc[linewidth=1.5pt](0,0){3.7}{120}{180}
\rput{0}(-0.2,3.5){$\alpha$}
\rput{0}(-0.2,1.5){$\beta$}
\rput{0}(2,1.8){$\gamma$}
\end{pspicture} }

\vskip-1.4cm

\centerline{\begin{pspicture}(-5,-1)(5,6)
\pscustom[fillstyle=solid, fillcolor=gray, linestyle=none]{%
  \psline(0,2.8284)(2.8284,2.8284)
  \psarc(0,0){4}{45}{90}
} %
\pscustom[fillstyle=solid, fillcolor=lightgray, linestyle=none]{%
  \psline(0,1)(3.8806,1)
  \psarc(0,0){4}{14.0362}{45}
  \psline(2.8284,2.8284)(0,2.8284)
} %
\pspolygon[fillstyle=solid, fillcolor=gray, linestyle=none](0,1)(3.8,1)(3.8,0.7071)(0,0.7071)
\pspolygon[fillstyle=solid, fillcolor=lightgray, linestyle=none](0,0.7071)(1.8708,0.7071)(0,0)
\pscustom[fillstyle=crosshatch, fillcolor=lightgray, linestyle=none]{%
  \psline(3.8,0)(2,0)
  \psarc(0,0){2}{0}{20.7048}
  \psline(1.8708,0.7071)(3.8,0.7071)
} %
\pscustom[fillstyle=solid, linestyle=none]{%
  \psline(4,0)(3.7,0)
  \psarc(0,0){3.7}{0}{20.7048}
  \psline(3.4610,1.3081)(3.7417,1.4142)
} %
\pscustom[fillstyle=solid, linestyle=none]{%
  \psarc(0,0){4}{0}{20.7048}
  \psline(3.7417,1.4142)(4,0)
} %
\pscustom[fillstyle=solid, linestyle=none]{%
  \psarc(0,0){1}{0}{180}
  \psline(0,1)(0,0)
  \psline(0,0)(1,0)
} %
\psarc(0,0){4}{0}{180}
\psarc(0,0){1}{0}{180}
\psline(-4,0)(-1,0)
\psline(1,0)(4,0)
\psline(0,1)(0,4)
\pscurve(0,2.8284)(0.2,2.5)(2,2)(3,1.4)(3.4,0)
\pscurve(-0,2.8284)(-0.2,2.5)(-2,2)(-3,1.4)(-3.4,0)
\psline[linestyle=dashed, linewidth=1.5pt, dash=10pt 2pt](-1,2.8284)(2.8284,2.8284)
\psline[linestyle=dashed, dash=2pt 4pt](2.8284,2.8284)(0.7071,0.7071)
\psline[linestyle=dashed, linewidth=1.5pt, dash=10pt 2pt](0.7071,0.7071)(3.9370,0.7071)
\psline[linewidth=1.5pt](0.9354,0.3536)(1.8708,0.7071)
\psarc[linewidth=1.5pt](0,0){2}{0}{20.7048}
\psline[linewidth=1.5pt](3.4610,1.3081)(3.7417,1.4142)
\psarc[linewidth=1.5pt](0,0){3.7}{0}{20.7048}
\psline[linewidth=1.5pt](-0.9354,0.3536)(-1.8708,0.7071)
\psarc[linewidth=1.5pt](0,0){2}{159.2952}{180}
\psline[linewidth=1.5pt](-3.4610,1.3081)(-3.7417,1.4142)
\psarc[linewidth=1.5pt](0,0){3.7}{159.2952}{180}
\rput{0}(-0.2,3.5){$\alpha$}
\rput{0}(-0.2,1.5){$\beta$}
\rput{0}(1.2,2){$\gamma$}
\end{pspicture} }

\vskip -.8cm
In each of these two pictures, a  tangential categorical open set is formed of the colored and hatched areas together with their symmetric relatively to $\alpha\cup \beta$. Denote it by $U$. It is deformed into a transverse set $A$ which is the union of the curves $\alpha$, $\beta$, $\gamma$. In the second picture $A$ lies inside $U$, showing that $A$ is an $A$-SND.
The tangential deformation of $U$ into $A$ is defined as follows:
\begin{itemize}
\item the dark gray zone of $U$ is deformed on the upper part $\alpha$ of the vertical line;
\item the clear gray zone is deformed on the lower part $\beta$ of the vertical line;
\item the hatched zone is deformed on the curve $\gamma$. 
\end{itemize}
\end{example}
\bibliographystyle{plain}
\bibliography{Tangcat}

\end{document}